\documentclass[12pt,a4paper,oneside,reqno,notitlepage]{amsart}

\usepackage[arrow,matrix,curve]{xy}\SilentMatrices
\def\xyma{\xymatrix@M.7em}
\usepackage{latexsym, amsbsy, amsmath, amsfonts, amssymb, amsthm, amscd, amscd}

\newtheorem{conj}{Conjecture}

\newtheorem{prop}{Proposition}
\newtheorem{theorem}{Theorem}
\newtheorem{lemma}{Lemma}

\begin{document}
\title{Intersection of subgroups in free groups and homotopy groups}

\author[Baues]{Hans-Joachim Baues}
\address{Max-Planck-Institut fur Mathematik, Vivatsgasse 7, 53111, Bonn, Germany}
\email{baues@mpim-bonn.mpg.de}

\author[Mikhailov]{Roman Mikhailov}
\address{Steklov Mathematical Institute, Gubkina 8, 119991 Moscow, Russia}
\email{rmikhailov@mail.ru}

\begin{abstract}
We show that the intersection of three subgroups in a free group
is related to the computation of the third homotopy group $\pi_3$.
This generalizes a result of Gutierrez-Ratcliffe who relate the
intersection of two subgroups with the computation of $\pi_2$. Let
$K$ be a two-dimensional CW-complex with subcomplexes
$K_1,K_2,K_3$ such that $K=K_1\cup K_2\cup K_3$ and $K_1\cap
K_2\cap K_3$ is the 1-skeleton $K^1$ of $K$. We construct a
natural homomorphism of $\pi_1(K)$-modules
$$
\pi_3(K)\to \frac{R_1\cap R_2\cap R_3}{[R_1,R_2\cap
R_3][R_2,R_3\cap R_1][R_3,R_1\cap R_2]},
$$
where $R_i=ker\{\pi_1(K^1)\to \pi_1(K_i)\},\ i=1,2,3$ and the
action of $\pi_1(K)=F/R_1R_2R_3$ on the right hand abelian group
is defined via conjugation in $F$. In certain cases, the defined
map is an isomorphism. Finally, we discuss certain applications of
the above map to group homology.
\end{abstract}

\thanks{The research of the second author was partially supported by
Russian Foundation Fond of Fundamental Research, grant N
05-01-00993 and Russian Presidental Grant MK-3166.} \maketitle

\section{Introduction}
\vspace{.5cm}

Simplicial homotopy theory makes it possible to translate certain
homotopy questions to the group-theoretical language. As a rule,
the group-theoretical problems appearing in this direction have a
difficult nature. Still combinatorial group theory is a crucial
tool in homotopy theory (see, for example, \cite{6A}, \cite{C}).
On the other hand, certain group-theoretical results may be
obtained by use of homotopy methods, like methods of simplicial
homotopy theory and the theory of derived functors (see, for
example, \cite{G}).

It is the purpose of this paper to combine the results of
Gutierrez-Ratcliffe and Wu which present certain links between
homotopical and group-theoretical structures. Let us recall them
first.\\

\noindent 1. (Exact sequence due to Guttierez-Ratcliffe,\
\cite{GR}) Let $K$ be a connected 2-dimensional CW-complex, and
$K_1$ and $K_2$ subcomplexes such that $K=K_1\cup K_2$ and
$K_1\cap K_2$ is the 1-skeleton $K^1$ of $K$, then there is an
exact sequence of $\pi_1(K)$-modules
\begin{equation}\label{exa1}
0\to i_1\pi_2(K_1)\oplus
i_2\pi_2(K_2)\buildrel{\alpha}\over{\rightarrow} \pi_2(K)\to
\frac{R\cap S}{[R,S]}\to 0,
\end{equation}
where $\alpha$ is induced by inclusion, $R$ is the kernel of
$\pi_1(K^1)\to\pi_1(K_1)$, $S$ is the kernel of
$\pi_1(K^1)\to\pi_1(K_1)$ and the action of $\pi_1(K)\simeq
\pi_1(K^1)/RS$ on $\frac{R\cap S}{[R,S]}$ is induced by
conjugation. The paper \cite{GR} contains another exact sequence
of the same nature. Let $G=\langle X\ |\ {\mathcal N}\rangle$ be a
group presentation with relation module $N/\gamma_2(N)$, where $N$
is the normal closure of the set $\mathcal N$ in the free group
$F(X)$. Let $K_r$ be a standard 2-complex for a presentation
$\langle X\ |\ r\rangle, r\in \mathcal N$, let $s_r$ be the root
of $r$ in the free group $F(X)$. There is an exact sequence of
$G$-modules:
\begin{equation}\label{exa2}
0\to \oplus_{r\in\mathcal
N}i_*\pi_2(K_r)\buildrel{\alpha}\over{\rightarrow} \pi_2(K)\to
\oplus_{r\in\mathcal N}\mathbb Z[G]/(s_r-1)\mathbb
Z[G]\buildrel{\gamma}\over{\rightarrow} N/\gamma_2(N)\to 0,
\end{equation}
where $\alpha$ is induced by inclusion and $\gamma$ maps
$1+(s_r-1)\mathbb Z[G]$ onto $r\gamma_2(N)$ for each $r\in\mathcal
N$. For a group-theoretical application of the above exact
sequences consider a two-relator presentation $\mathcal P=\langle
X\ |\ r_1,r_2\rangle$ of a group $F(X)/\langle
r_1,r_2\rangle^{F(X)}$. As a consequence of the exact sequences
(\ref{exa1}) and (\ref{exa2}), we have that the quotient
$\frac{\langle r\rangle^{F(X)}\cap \langle
r_2\rangle^{F(X)}}{[\langle r_1\rangle^{F(X)},\,\langle
r_2\rangle^{F(X)}]}$ is a subgroup of the quotient $\mathbb
Z[G]/(s_{r_1}-1)\mathbb Z[G]\oplus \mathbb Z[G]/(s_{r_2}-1)\mathbb
Z[G]$, i.e. it is a free Abelian. Hence, we proved the following
generalization of the theorem due to Hartley-Kuzmin \cite{HK}:
{\it let $F$ be a free group, $r_1$ and $r_2$ words in $F$,
$R_i=\langle r_i\rangle^F,\ i=1,2$, then the group $\frac{R_1\cap
R_2}{[R_1,\,R_2]}$ is a free Abelian group.}
\\

\noindent 2. (The presentation of homotopy groups of the 2-sphere
due to Wu) It is one of the deep problems of algebraic topology to
compute homotopy groups $\pi_n(S^2)$. Note that the 2-sphere
presents the most 'unstable' case from the point of view of
homotopy theory. The developed methods of Adams-type spectral
sequences usually do not work in this case. In low degrees one has
(see \cite{To}):
\begin{center}
\begin{tabular}{cccccccccc}
 $n$ & \vline & 2 & 3 & 4 & 5 & 6 & 7 & 8 & 9 \\ \hline  $\pi_n(S^2)$ & \vline &
 $\mathbb Z$ & $\mathbb Z$ & $\mathbb Z_2$ & $\mathbb Z_2$ & $\mathbb
 Z_4\oplus \mathbb Z_3$ & $\mathbb Z_2$ & $\mathbb Z_2$ & $\mathbb
 Z_3$
\end{tabular}
\end{center}
The structure of $\pi_n(S^2)$ is known up to some stage ($n\approx
30$), mostly due to Toda and his students. The general structure
of $\pi_n(S^2)$ is unclear and mysterious.

Recall the description of homotopy groups of the 2-sphere due to
Wu \cite{Wu}. Let $F[S^1]$ be Milnor's $F[K]$-construction applied
to the simplicial circle $S^1$. This is the free simplicial group
with $F[S^1]_n$ a free group of rank $n\geq 1$ with generators
$x_0,\dots, x_{n-1}$. Changing the basis of $F[S^1]_n$ in the
following way: $y_i=x_ix_{i+1}^{-1},\ y_{n-1}=x_{n-1}$, we get
another basis $\{y_0,\dots,y_{n-1}\}$ in which the simplicial maps
can be written easier. A combinatorial group-theoretical argument
then gives the following description of the $n$-th homotopy group
of the loop space $\Omega\Sigma S^1$, which is isomorphic to the
homotopy group of $\pi_{n+1}(S^2)$ (see \cite{Wu} for explicit
computations):
\begin{equation}\label{wujie}
\pi_{n+1}(S^2)\cong \frac{\langle y_{-1}\rangle^F\cap \langle
y_0\rangle^F \cap\dots\cap \langle
y_{n-1}\rangle^F}{[[y_{-1},y_0,\dots, y_{n-1}]]},\ n\geq 1
\end{equation}
where $F$ is a free group with generators $y_0,\dots, y_{n-1},$
$y_{-1}=(y_0\dots y_{n-1})^{-1},$ the group $[[y_{-1},y_0,\dots,
y_{n-1}]]$ is the normal closure in $F$ of the set of left-ordered
commutators \begin{equation}\label{comm}
[z_1^{\varepsilon_1},\dots, z_{t}^{\varepsilon_t}]\end{equation}
with the properties that $\varepsilon_i=\pm 1$,
$z_i\in\{y_{-1},\dots, y_{n-1}\}$ and all elements in
$\{y_{-1},\dots, y_{n-1}\}$ appear at least once in the sequence
of elements $z_i$ in (\ref{comm}).\\

The main idea of our approach can be formulated as the following
{\it conjectural observation}: {\it the nature of the presentation
(\ref{wujie}) comes from the fact that the 2-sphere is
homotopically equivalent to the standard 2-complex, constructed
from the presentation $\langle y_0,\dots, y_{n-1}\ |\ y_0,\dots,
y_{n-1},y_{n-1}^{-1}\dots y_0^{-1}\rangle$.} Given a free group
$F$ and normal subgroups $(n\geq 2)$
$$
R_1,\dots, R_n\unlhd F,
$$
denote the quotient group
$$
I_n(F,R_1,\dots, R_n):=\frac{R_1\cap\dots \cap R_n}{\prod_{I\cup
J=\{1,\dots, n\},\ I\cap J=\emptyset}[\bigcap_{i\in I}R_i,
\bigcap_{j\in J}R_j]}.
$$
Here $\bigcap$ denotes the intersection of subgroups in the free
group $F$ and $\prod$ is the product of commutator subgroups as
indicated. In fact, the abelian group $I_n$ has the natural
structure of an $F/R_1\dots R_n$-module, with the group action
defined via conjugation in $F$.

The computation of the abelian group $I_n$ is highly non-trivial.
In fact, in the special case $F=\langle x_1,\dots, x_{n-1}\rangle,
R_i=\langle x_i\rangle^F,\ i=1,\dots, n-1, R_n=\langle x_1\dots
x_{n-1}\rangle^F$ a standard commutator calculus argument, given
essentially in Corollary 3.5 of \cite{Wu} shows that
$$
[[y_{-1},y_0,\dots, y_{n-1}]]=\prod_{I\cup J=\{1,\dots, n+1\},\
I\cap J=\emptyset}[\bigcap_{i\in I}R_i,\bigcap_{j\in J}R_j],
$$
and hence we have the following isomorphism
$$
I_n(F,R_1,\dots, R_n)=\pi_n(S^2).
$$

On the other hand, for $n=2$, one has a general description of the
$F/R_1R_2$-module $I_2(F,R_1,R_2)=\frac{R_1\cap R_2}{[R_1,\,R_2]}$
in terms of homotopy groups of certain spaces given by the
sequence (\ref{exa1}) due to Gutierrez-Ratcliffe. For the
generalization of the Gutierrez-Ratcliffe's approach to the higher
dimensional homotopy groups  consider a connected 2-dimensional
CW-complex $K$ with subcomplexes $$ K_1,\dots, K_n\subset K,
$$
for which $K_1\cup\dots\cup K_n=K$ and $K_1\cap\dots \cap K_n$ is
the 1-skeleton $K^1$ of $K$, with $F=\pi_1(K^1)$ and
$$
R_i=ker\{\pi_1(K^1)\to \pi_1(K_i)\},\ i=1,\dots,n.
$$

We conjecture that each element $\alpha\in \pi_n(S^2)$ determines
a natural function $(n\geq 2)$
$$
\alpha_*: \pi_2(K)/(i_1\pi_2(K_1)+\dots +i_n\pi_2(K_n))\to
I_n(F,R_1,\dots, R_n).
$$
In general, $\alpha_*$ is not a homomorphism of abelian groups.\\

 \noindent{\bf Proposition.} {\it Let $n=2$. If $\alpha$ is a
generator of $\pi_2(S^2)=\mathbb Z$, then $\alpha_*$ exists and is
given by the map $\pi_2(K)\to I_2(F,R_1,R_2)$ of
Gutierrez-Ratcliffe \cite{GR}.}
\\ \\
Moreover, as a main result of this paper we prove the following\\
\\  \noindent{\bf Theorem.} {\it Let $n=3$. If $\alpha\in \pi_3(S^2)$ is
a generator, then there is a well-defined function $\alpha_*$
which is a quadratic map inducing a natural homomorphism of
$\pi_1(K)$-modules
$$
\alpha_\#:\pi_3(K)\to I_3(F,R_1,R_2,R_3).
$$
For the example of Wu, one has $K=S^2$ and in this case
$\alpha_\#$ is an isomorphism.}\\

In the construction of the above homomorphism of
$\pi_1(K)$-modules we essentially use the fact that
$\pi_3(K)=\Gamma\pi_2(K),$ where $\Gamma$ is Whitehead's universal
quadratic functor.

One can interpret certain elements of $\pi_*(S^2)$ as the elements
of certain free groups in Milnor's $F[S^1]$-construction. For
example, the element in $F(y_0,y_1,y_2)$ corresponding to the
generator of $\pi_4(S^2)$ in (\ref{wujie}) is
$$
[[y_0,y_1],[y_0,y_1y_2]].
$$
The element in $F(y_0,y_1,y_2,y_3)$, corresponding to the
generator of $\pi_5(S^2)$ is
$$
[[[y_0,y_1],[y_0,y_1y_2]],[[y_0,y_1],[y_0,y_1y_2y_3]]].
$$
This follows from the result of Wu \cite{Wu}. In Section 3 we
shall formulate a conjecture, which, when applied to these
elements, leads to non-trivial commutator problems in free groups
that we are not ready to solve. This difficulty is the reason why
we consider in this paper only the case of three subcomplexes.

Finally, we use the main construction of the paper for an
arbitrary free group $F$ and its normal subgroups $R_1,R_2,R_3$,
to define the following natural map of abelian groups:
$$
H_4(G)\to \frac{R_1\cap R_2\cap R_3}{[R_1,R_2\cap R_3][R_2,R_3\cap
R_1][R_3,R_1\cap R_2][F,R_1\cap R_2\cap R_3]}.
$$
This map is related to the Brown-Ellis construction (see
\cite{BE}), however, the methods used in the current paper differ
from those used in \cite{BE}. The relation with Brown-Ellis
construction is given in \cite{EM}.\vspace{.5cm}
\section{The category $\mathcal K_n$}
\vspace{.5cm} For $n\geq 2,$ denote by $\mathcal K_n$ the category
with objects $\bar K=(K, K_1,\dots, K_n).$ Here $K$ is a
two-dimensional CW-complex, $K_i$ is a subcomplex of $K$,
$i=1,\dots, n$, such that $K=K_1\cup\dots \cup K_n,$  and
$K^1=K_1\cap \dots\cap K_n.$ A morphism in $Hom_{\mathcal
K_n}(\bar K,\bar L)$ for $\bar K,\bar L\in \mathcal K_n$ is a map
$$
f: K^1\to L^1
$$
between 1-skeletons of $K$ and $L$, such that $f$ can be extended
to a map $\bar f: K\to L$, with the property $\bar f(K_i)\subseteq
L_i,\ i=1,\dots, n$.

Denote by $\mathcal R_n$ $(n\geq 2)$ the category with objects
$(F,R_1,\dots, R_n),$ where $F$ is a free group and $R_i$ is a
normal subgroup in $F$. A morphism in $\mathcal R_n$ between two
objects $(F, R_1,\dots, R_n)$ and $(F', R_1',\dots, R_n')$ is a
group homomorphism $g: F\to F'$ such that $g(R_i)\subseteq R_i',\
i=1,\dots, n$. This category was also considered in \cite{DIP}.

There is a natural functor between these two categories,
$$\mathcal F_n: \mathcal K_n\to \mathcal R_n,$$ defined by setting
$$
\mathcal F_n: (K,K_1,\dots, K_n)\mapsto (\pi_1(K^1), R_1,\dots,
R_n),
$$
where $R_i=\ker\{\pi_1(K^1)\to \pi_1(K_i)\}.$

For $n\geq 2$, define the functor
$$
I_n: \mathcal R_n\to \mathcal Ab,
$$
where $\mathcal Ab$ is the category of abelian groups, by setting
$$
I_n: \bar R=(F,R_1,\dots, R_n)\mapsto I_n(\bar
R):=\frac{R_1\cap\dots \cap R_n}{\prod_{I\cup J=\{1,\dots, n\},\
I\cap J=\emptyset}[\bigcap_{i\in I}R_i, \bigcap_{j\in J}R_j]}.
$$
Clearly, for any $\bar R\in \mathcal R_n$, the abelian group
$I_n(\bar R)$ has a natural structure of $F/R_1\dots R_n$-module,
where the group action viewed via conjugation in $F$.

\vspace{.5cm}
\section{The surjection $q$ and the conjecture on $\alpha_*$}
\vspace{.5cm}
\subsection{}
Consider the two-dimensional sphere $S^2$ as the standard
two-complex constructed from the following presentation of the
trivial group:
\begin{equation}\label{oe}
\langle x_1,\dots, x_{n-1}\ |\ x_1,\dots, x_{n-1},\
x_{n-1}^{-1}\cdots x_1^{-1}\rangle.
\end{equation}
This presentation defines an element $\bar S_n$ from $\mathcal
K_n$:
\begin{equation}\label{sbar}
\bar S_n=(S^2, L_1,\dots, L_n),
\end{equation}
with $L_i=\vee_{i=1}^{n-1}S^1\cup e_i,$  where $e_i$ is the 2-cell
corresponding to the relation word $x_i,\ i=1,\dots, n-1,$ $e_n$
is the 2-cell corresponding to the relation word
$x_{n-1}^{-1}\cdots x_1^{-1}$.

In this section we show the following result.

\begin{prop}
For an object $\bar S_n$ in $\mathcal K_n$ associated to Wu's
example in $\mathcal R_n$ there is a surjection
$$
q: Hom_{\mathcal K_n}(\bar S_n,\bar K)\twoheadrightarrow
\pi_2(K)/(i_1(K_1)+\dots i_n(K_n)),
$$
which is natural in $\bar K\in \mathcal K_n$.
\end{prop}

For $\alpha\in \pi_n(S^2)=I_n\mathcal F_n(\bar S_n)$ we thus
obtain the following diagram

$$
\xyma{ &Hom_{\mathcal K_n}(\bar S_n,\bar K)\ar@{->}[dd]^{\alpha^*}
\ar@{->>}[rr]&{}^q&\pi_2(K)/(i_1\pi_2(K_1)+\dots +i_n\pi_2(K_n))
\ar@{-->}[ddll]_{\alpha_*} \\
\\ & I_n(\mathcal F_n(\bar K)) && \\}
$$
where $\alpha^*(f)=f_*(\alpha)$.

\begin{conj}
For each $\alpha\in \pi_n(S^2)$ there exists a function $\alpha_*$
for which the diagram commutes. Hence $\alpha_*$ is well defined
and natural provided $q(f)=q(g)$ implies
$\alpha^*(f)=\alpha^*(g)$.
\end{conj}

\subsection{}
Recall that for a given two-dimensional complex $K$, the free
crossed module $$\partial: \pi_2(K,K^1)\to \pi_1(K^1)$$ can be
defined as follows. The group $\pi_2(K, K^1)$ is generated by the
set $$\{e_{\alpha}^w\ |\ \alpha\ \text{is a 2-cell in}\ K,\ w\in
\pi_1(K^1)\}$$ with the set of relations
\begin{equation}\label{relset}\{e_{\alpha}^ve_{\beta}^we_{\alpha}^{-v}e_\beta^{-u},\
u=vr_\alpha v^{-1}w\},\end{equation} where $r_\alpha\in
\pi_1(K^1)$ is the attaching element representing $e_\alpha$ (see,
for example, \cite{GH}). The homomorphism $\partial$ is defined by
setting $\partial: e_{\alpha}^w\mapsto r_\alpha^w$. Hence every
element from $ker(\partial)=\pi_2(K)$ can be represented by an
element $e_{\alpha_1}^{\pm w_1}\dots e_{\alpha_m}^{\pm w_m},$ such
that $r_{\alpha_1}^{\pm w_1}\dots r_{\alpha_m}^{\pm w_m}$ is
trivial in $\pi_1(K^1)$.

Let $f\in Hom_{\mathcal K_n}(\bar S_n, \bar K)$. It means that
there exists a homomorphism between two free groups $f:
F_{n-1}:=F(x_1,\dots, x_n)\to \pi_1(K^1)$ such that
\begin{equation}\label{c1}f(x_i)\in \ker\{\pi_1(K^1)\to
\pi_i(K_i)\},\ i=1,\dots, n-1\end{equation} and $f$ can be
extended to a homomorphism between two crossed modules:
\begin{equation}\label{dia}
\begin{CD}
\pi_2(S^2, \vee_{i=1}^{n-1}S^1)@>\partial_1>> F_{n-1}\\
@V{f'}VV @V{f}VV\\
\pi_2(K,K^1) @>\partial_2>> \pi_1(K^1)
\end{CD}
\end{equation}
For a given group homomorphism $f: F_{n-1}\to \pi_1(K^1)$ with the
property (\ref{c1}), the necessary and sufficient condition for
the existence of the extension (\ref{dia}) is the condition
$$
f(x_1\cdots x_n)\subseteq R_n:=\ker\{\pi_1(K^1)\to \pi_1(K)\}.
$$

For $\bar K=(K,K_1,\dots, K_n)\in \mathcal K_n,$ we now define the
canonical (forgetful) map
$$
q: Hom_{\mathcal K_n}(\bar S_n, \bar K)\to
\pi_2(K)/(i_1\pi_2(K_1)+\dots i_n\pi_2(K_n)),
$$
which carries a morphism $\bar S^2\to \bar K$ to the underlying
map $S^2\to K$. Here the natural maps $i_j: \pi_2(K_j)\to K$ are
induced by inclusions $K_j\to K$. Using the language of crossed
modules, we can describe the map $q$ as follows. Denote by
$\{s_1,\dots, s_{n}\}$ the set of 2-cells in $S^2$ viewed as the
standard two-complex for the group presentation (\ref{oe}). The
map $f'$ defines elements $f'(s_\alpha)\in \pi_2(K,K^1).$ Observe
that $\partial_1(s_1\dots s_n)=1$ and the element $s_1\dots s_n$
presents the generator of $\pi_2(S^2)$. Since the diagram
(\ref{dia}) is commutative, $\partial_2(f'(s_1)\dots f'(s_n))=1$
and the element $f'(s_1)\dots f'(s_n)$ represents certain element
from $ker(\partial_2)=\pi_2(K)$, which is exactly $q(f)$. Let us
show that this map does not depend on an extension (\ref{dia}).
Suppose we have another extension of the homomorphism $f$:
\begin{equation}\label{dia2}
\begin{CD}
\pi_2(S^2, \vee_{i=1}^{n-1}S^1)@>\partial_1>> F_{n-1}\\
@V{f''}VV @V{f}VV\\
\pi_2(K,K^1) @>\partial_2>> \pi_1(K^1)
\end{CD}
\end{equation}
with $f''(s_j)\neq f'(s_j)$ at least for one $j$ $(1\leq j\leq
n)$. It follows that $\partial_2(f'(s_j)f''(s_j)^{-1})=1,$ hence
$$
f'(s_j)f''(s_j)^{-1}\in im\{i_j: \pi_2(K_j)\to \pi_2(K)\}.
$$
Therefore, the images of elements $f'(s_1\dots s_n)$ and
$f''(s_1\dots s_n)$ are equal in the quotient
$\pi_2(K)/(i_1\pi_2(K_1)+\dots i_n\pi_2(K_n))$ and the map $q$ is
well-defined.
\begin{lemma}
The map $q$ is surjective.
\end{lemma}
\begin{proof}
Consider the diagram (\ref{dia2}). Now let $c=e_{\alpha_1}^{\pm
w_1}\dots e_{\alpha_m}^{\pm w_m}$ be an arbitrary element from
$ker(\partial_2)$. Let us enumerate all cells of $K$ in the
following order: $e_{1,\alpha},\dots, e_{n,\alpha}$ with
$e_{i,\alpha}\in K_i,\ i=1,\dots, m$. Clearly, the set of
relations (\ref{relset}) in $\pi_2(K,K^1)$ gives a possibility to
present the element $c$ in the form
$$
c=\prod_{*} e_{1,*}^{\pm w_{1,*}}\dots \prod_* e_{n,*}^{\pm
w_{n,*}}
$$
with some $w_{i,*}\in \pi_1(K^1)$. We define the map $f:
F_{n-1}\to \pi_1(K^1)$ by setting $f(x_i)=\prod_* r_{i,*}^{\pm
w_{n,*}}$. We can extend it to $f':
\pi_2(S^2,\vee_{i=1}^{n-1}S^1)\to \pi_2(K,K^1)$ by
$f'(s_i)=\prod_* r_{i,*}^{\pm w_{n,*}}$. This is correct, since $$
\partial_1(f'(s_n))=\partial_2(f'(s_1)\dots f'(s_{n-1}))^{-1}=f(\partial_1(s_1\dots
s_{n-1})^{-1}).
$$
The homotopy class corresponding to the element $c\in
\pi_2(K,K^1)$ coincides with $q(f)$ and the surjectivity of $q$ is
proved.
\end{proof}
\vspace{.5cm}
\section{Proof of the conjecture for $n=2$ and $n=3$}
\vspace{.5cm}
\subsection{}
There are different ways of description of elements from $\pi_2$
for a standard complex of a given group presentation, for example,
pictures, kernels of Jacobian maps, defined via Fox calculus etc.
We describe the map $q$ in the conjecture by use of identity
sequences which represent elements in $\pi_2(K)$, see \cite{Pr}.
The material about identity sequences we recall here is
well-known. Nevertheless, we give it in detail since it will be
the basic technical devise in our proofs.

Let $F$ be a free group with basis $X$ and $\mathcal R$ a certain
set of words in $F$. Consider the group presentation
\begin{equation}\label{pres} \mathcal P=\langle X\ |\ \mathcal R\rangle
\end{equation}
$c_i,\ i=1,\dots, m$ are words in $F$, which are conjugates of
elements from $\mathcal R$, i.e. $c_i=t_i^{\pm w_i},\ t_i\in
\mathcal R,\ w_i\in F$. The sequence
\begin{equation}\label{seqident}
c=(c_1,\dots, c_m)
\end{equation}
is called an {\it identity sequence} if the product $c_1\dots c_m$
is the identity in $F$. For a given identity sequence
(\ref{seqident}), define its inverse:
$$
c^{-1}=(c_m^{-1},\dots, c_1^{-1}).
$$
For a given element $w\in F$, the conjugate $c^w$ is the sequence:
$$
c^w=(c_1^w,\dots, c_m^w),
$$
which clearly is again an identity sequence. Define the following
operation in the class of identity sequences, called {\it Peiffer
operations}:\\
(i) replace each $w_i$ by any word equal to it in $F$;\\
(ii) delete two consecutive terms in the sequence if one is equal
identically to the inverse of the other;\\
(iii) add two consecutive terms in the sequence if one is equal
identically to the inverse of the other;\\
(iv) replace two consecutive terms $c_i, c_{i+1}$ by terms
$c_{i+1},c_{i+1}^{-1}c_ic_{i+1}$;\\
(v) replace two consecutive terms $c_i, c_{i+1}$ by terms
$c_ic_{i+1}c_i^{-1}, c_i$.

Two identity sequences are called {\it equivalent} if one can be
obtained from the other by a finite number of Peiffer operations.
This defines an equivalence relation in the class of identity
sequences. The set of equivalence classes of identity sequences
for a given group presentation (\ref{pres}) denote by $E_{\mathcal
P}$. The set $E_{\mathcal P}$ can be viewed as a group, with a
binary operation defined as a class of justaposition of two
sequences: for identity sequences $c_1,c_2$ and their equivalence
classes $\langle c_1\rangle, \langle c_2\rangle \in E_{\mathcal
P}$, $\langle c_1\rangle +\langle c_2\rangle =\langle
c_1c_2\rangle$. The inverse element of the class $\langle
c\rangle$ is $\langle c^{-1}\rangle$ and the identity in
$E_{\mathcal P}$ is the empty sequence. It is easy to see that
$E_{\mathcal P}$ is Abelian. For two identity sequences
$c=(c_1,\dots, c_m)$ and $d=(d_1,\dots, d_k)$, we have
$$
\langle cd\rangle=\langle (c_1,\dots, c_m, d_1,\dots,
d_k)\rangle=\langle (d_1,\dots, d_k, c_1^{d_1\dots d_m},\dots,
c_m^{d_1\dots d_m})\rangle
$$
by the relation (iv). Since $d_1\dots d_m=1$ in $F$, we have
$$
\langle cd\rangle=\langle (d_1,\dots, d_k, c_1,\dots, c_m)\rangle
=\langle dc\rangle.
$$
Furthermore, $E_{\mathcal P}$ is a $F$-module, where the action is
given by
$$
\langle c\rangle \circ f=\langle c^f\rangle,\ f\in F.
$$
It is easy to show that
$$
\langle c\rangle \circ r=\langle c\rangle,\ r\in R,
$$
i.e. the subgroup $R$ acts trivially at $E_{\mathcal P}$. To see
this, let $r=r_1^{\pm v_1}\dots r_k^{\pm v_k},\ r_i\in \mathcal
R,\ v_i\in F$. For any identity sequence $c=(c_1,\dots, c_m)$, by
(ii), (iii), (iv),
\begin{multline*}
\langle (c_1,\dots, c_m)\rangle= \langle (c_1,\dots, c_m, r_1^{\pm
v_1},\dots, r_k^{\pm v_k}, r_k^{\mp v_k},\dots, r_1^{\mp v_1})\rangle=\\
\langle (r_1^{\pm v_1},\dots, r_k^{\pm v_k}, c_1^r,\dots, c_m^r,
r_k^{\mp v_k},\dots, r_1^{\mp v_1})\rangle=\langle (c_1^r,\dots,
c_m^r)\rangle.
\end{multline*}
Thus $E_{\mathcal P}$ can be viewed as a $G$-module. It is not
hard to show that for a given presentation $\mathcal P$, the
second homotopy module $\pi_2(K_\mathcal P)$ is isomorphic to the
identity sequence module $E_{\mathcal P}$ (see, for example,
\cite{Pr}).

\subsection{} For a given $\bar K$ choose the elements
$e_{i,\alpha}\in \pi_2(K,K^1),\ i=1,\dots, n,\alpha\in A$ which
represent the corresponding two-dimensional cells in $K_i, \
i=1,\dots, n$ with the natural property
$$
\partial(e_{i,\alpha})\in R_i,
$$
where $R_i=\ker\{\pi_1(K^1)\to \pi_1(K_i)\},\ i=1,\dots n,$ and
the normal closure of the set $\{\partial(e_{i,\alpha})\ |\
\alpha\in A\}$ in $\pi_1(K^1)$ is equal to $R_i$. Clearly, $K$ is
homotopically equivalent to a wedge
$$
K\simeq \bigvee_{j\in J}S^2\vee K_{\mathcal P},
$$
where $K_{\mathcal P}$ is the standard two-complex constructed
from the group presentation
$$
\langle X\ |\ \partial(e_{i,\alpha}),\ i=1,\dots, n,\ \alpha\in
A\rangle,
$$
with $X$ being a basis of $\pi_1(K^1)$. We have the following
natural isomorphism of $\pi_1(K)$-modules:
$$
\pi_2(K)/(i_1\pi_2(K_1)+\dots+i_n\pi_2(K_n))\simeq
\pi_2(K_{\mathcal P})/(i_1\pi_2(K_{\mathcal
P_1})+\dots+i_n\pi_2(K_{\mathcal P_n})),
$$
where $\mathcal P_i$ is the following presentation of the group
$\pi_1(K_i)$:
$$
\langle X\ |\ \partial(e_{i,\alpha}),\ \alpha\in A\rangle
$$
for $i=1,\dots, n$.

Let $f,g\in Hom_{\mathcal K_n}(\bar S_n, \bar K)$. We can present
\begin{align*}& f(x_i)={r_1^{(i)}}^{\pm w_{1,i}}\dots
{r_{k_i}^{(i)}}^{\pm w_{k_i,i}},\ i=1,\dots, n-1,\\
& f(x_1\cdots x_n)={r_1^{(n)}}^{\pm w_{1,n}}\dots
{r_{k_n}^{(n)}}^{\pm w_{k_n,n}}
\end{align*}
for some $r_{j}^{(i)}\in \{\partial(e_{i,\alpha}),\ \alpha\in A\}$
and $w_{j,i}\in \pi_1(K^1)$. Analogically for $g\in Hom_{\mathcal
K_n}(\bar S_n, \bar K)$:
\begin{align*}& g(x_i)={r_1'^{(i)}}^{\pm w_{1,i}'}\dots
{r_{k_i'}'^{(i)}}^{\pm w_{k_i',i}'},\ i=1,\dots, n-1,\\
& g(x_1\cdots x_n)={r_1'^{(n)}}^{\pm w_{1,n}'}\dots {
r_{k_n'}'^{(n')}}^{\pm w_{k_n',n}'}
\end{align*}

The following Lemma follows directly from the definition of the
map $q$ and the above description of the second homotopy module
for the standard complex in terms of identity sequences.

\begin{lemma}
Using the above notation, $q(f)=q(g)$ if and only if the identity
sequence
$$
({r_1^{(1)}}^{\pm w_{1,1}},\dots, {r_{k_n}^{(n)}}^{\pm w_{k_n,n}},
{ r_{k_n'}'^{(n')}}^{\mp w_{k_n',n}'},\dots, \dots,
{r_{1}'^{(1)}}^{\mp w_{1,1}'})
$$
is equivalent to an identity sequence of the form
$$
(s_1^{(1)},\dots, s_{l_1}^{(1)}, \dots, s_1^{(n)},\dots
,s_{l_n}^{(n)})
$$
with $s_j^{(i)}\in \{\partial(e_{i,\alpha})^{\pm w},\
w\in\pi_1(K^1)\}$ such that $s_1^{(i)}\dots s_{l_i}^{(i)}$ is
trivial in $\pi_1(K^1)$ for every $i=1,\dots, n$.
\end{lemma}

Let $(K,K_1,K_2)\in \mathcal K_2$. The $\pi_1(K)$-module
$\pi_2(K)/(i_1\pi_2(K_1)+i_2\pi_2(K_2))$ can be identified to the
module of the identity sequences of the type
\begin{equation}\label{sequ}(c_1,\dots, c_m),\ c_j\in \{\partial(c_{i,\alpha})^w,\
w\in \pi_1(K^1),\ \alpha\in A,\ i=1,2\}\end{equation} modulo the
sequences of the form $(c_1,\dots, c_{m_1},c_{m_1+1},\dots, c_m)$
with $$c_1,\dots, c_{m_1}\in \{\partial(c_{1,\alpha})^w,\ w\in
\pi_1(K^1)\}, c_{m_1+1},\dots, c_{m}\in
\{\partial(c_{2,\alpha})^w,\ w\in \pi_1(K^1)\}$$ and
$$
c_1\dots c_{m_1}=c_{m_1+1}\dots c_m=1
$$
in $\pi_1(K^1)$.

Every identity sequence (\ref{sequ}) with the help of Peiffer
operations of the type (iv) can be reduced to the sequence of the
form $(c_1,\dots, c_{m_1},c_{m_1+1},\dots, c_m)$ with $c_1,\dots,
c_{m_1}\in \{\partial(c_{1,\alpha})^w,\ w\in \pi_1(K^1)\},$
$c_{m_1+1},\dots, c_{m}\in \{\partial(c_{2,\alpha})^w,\ w\in
\pi_1(K^1)\}.$

\subsection{}
For the most elementary case $n=2$ we view the 2-sphere $S^2$ as a
standard complex constructed from the group presentation
$$
\langle x\ |\ x, x^{-1}\rangle.
$$
Clearly then
$$
I_2(\bar S^2)=\frac{\langle x\rangle\cap \langle
x^{-1}\rangle}{[\langle x\rangle, \langle x\rangle]}\simeq \mathbb
Z
$$
with $x$ a generator of this infinite cyclic group. For the
generator $x\in \pi_2(S^2)$, the map
$$
\Lambda_x: \pi_2(K)/(i_1\pi_2(K_1)+i_2\pi_2(K_2))\to \frac{R_1\cap
R_2}{[R_1,R_2]}
$$
is given in the above notation by
$$
\Lambda_x: (c_1,\dots, c_{m_1},c_{m_1+1},\dots, c_m)\mapsto
c_1\cdots c_{m_1}.[R_1,R_2].
$$
First observe that $\Lambda_x$ is the homomorphism of
$\pi_1(K)=\pi_1(K^1)/R_1R_2$-modules. Secondly, $\Lambda_x$
clearly is an epimorphism. The fact that $\Lambda_x$ is a
monomorphism is not difficult (see Theorem 1.3 \cite{Pr} for the
complete proof). Hence we have the following exact sequence of
$\pi_1(K)$-modules due to Gutierrez and Ratcliffe \cite{GR}:
\begin{equation}\label{grs}
0\to i_1\pi_2(K_1)+
i_2\pi_2(K_2)\buildrel{\alpha}\over{\rightarrow} \pi_2(K)\to
\frac{R_1\cap R_2}{[R_1,R_2]}\to 0.
\end{equation}

\begin{theorem}
Conjecture 1 is true for $n=3$.
\end{theorem}
\begin{proof}
In this case we view $S^2$ as the standard complex constructed for
the group presentation
$$
\langle x_1,x_2\ |\ x_1,x_2,x_2^{-1}x_1^{-1}\rangle
$$
with
$$
I_3(\mathcal F_3(\bar S^2))=I_3(F(x_1,x_2),\langle
x_1\rangle^{F(x_1,x_2)}, \langle x_2\rangle^{F(x_1,x_2)}, \langle
x_2^{-1}x_1^{-1}\rangle^{F(x_1,x_2)})\simeq \mathbb Z
$$
with a generator given by the commutator $[x_1,x_2]$.

\subsection{}
Let $\bar K=(K,K_1,K_2, K_3)\in \mathcal K_3$. Denote $F=\pi_1(K^1)$. Denote the sets of words in $F$:
$
\mathcal R_i=\{\partial(e_{i,\alpha},\ \alpha\in A\},\ i=1,2,3.
$
By $\mathcal R_i^F$ we mean the set $\{r^w,\ r\in \mathcal R_i,\ w\in F\}.$
The $\pi_1(K)$-module
$\pi_2(K)/(i_1\pi_2(K_1)+i_2\pi_2(K_2)+i_3\pi_2(K_3))$ can be
identified with the module of the identity sequences
\begin{equation}\label{sequ1}c=(c_1,\dots, c_m),\ c_j\in \mathcal R_1^F\cup \mathcal R_2^F\cup \mathcal R_3^F\end{equation}
modulo the sequences of the type
\begin{equation}\label{sequ3}(c_1,\dots, c_{m_1},c_{m_1+1},\dots, c_{m_2},c_{m_2+1},\dots,
c_m)\end{equation}
with
$c_1,\dots, c_{m_1}\in \mathcal R_1^F,$ $c_{m_1+1},\dots,
c_{m}\in \mathcal R_2^F,$
$c_{m_2+1},\dots, c_m\in \mathcal R_3^F$ and
\begin{equation}\label{ccond}
c_1\dots c_{m_1}=c_{m_1+1}\dots c_{m_2}=c_{m_2+1}\dots c_m=1,
\end{equation}
in $F$.

Divide the sequence
(\ref{sequ1}) into the three ordered subsequences
\begin{equation}\label{div}
(c_{r_1},\dots, c_{r_l}),\ \ (c_{s_1},\dots, c_{s_k}),\ \
(c_{t_1},\dots, c_{t_n}), \end{equation} where $c_{r_i}\in
\mathcal R_1^F,\ i=1,\dots, l$,
$c_{s_i}\in \mathcal R_2^F,\ i=1,\dots,
k$, $c_{t_i}\in\mathcal R_3^F,\ i=1,\dots, h$ and $$ r_1<r_2<\dots<
r_l,\ \ s_1<s_2<\dots<s_k,\ \ t_1<t_2<\dots<t_h,
$$
$$
\{r_1,\dots, r_l\}\cup\{s_1,\dots, s_k\}\cup\{t_1,\dots,
t_h\}=\{1,\dots, m\}.
$$
Denote
\begin{align*}
& \bar c_i=c_{r_i},\ i=1,\dots, l,\\
& \bar c_{l+i}=c_{s_i}^{\prod_{r_j>s_i}c_{r_j}},\ i=1,\dots, k,\\
& \bar
c_{l+k+i}=c_{t_1}^{(\prod_{r_z>t_i}c_{r_z})\prod_{s_j>t_1}c_{s_j}^{\prod_{r_z>s_j}c_{r_z}}},\
i=1,\dots, h.
\end{align*}
Clearly, $$\bar c_1,\dots, \bar c_l\in R_1,\ \bar c_{l+1},\dots,
\bar c_{l+k}\in R_2,\ \bar c_{l+k+1},\dots, \bar c_{l+k+h}\in R_3$$
and the sequence
\begin{equation}\label{wew}
(\bar c_1,\dots, \bar c_{l+k+h})
\end{equation}
is made of the sequence (\ref{sequ1}), applying the Peiffer operations of
type (iv). At the first step we replace all terms $c_{r_i}$ to the
left side of the sequence. At the second step we replace all terms
$c_{s_i}$ between elements $c_{r_i}$-s and $c_{t_i}$-s and get the
sequence (\ref{wew}). Denote
\begin{align*} &  r_c:=\bar c_{1}\dots \bar c_{l}\in R_1,\\
& s_c:=\bar c_{l+1}\dots \bar c_{l+k}\in R_2,\\
& t_c:=\bar c_{l+k+1}\dots \bar c_{l+k+h}\in R_3.
\end{align*}
In these notations, for the generator $x:=[x_1,x_2]$ of $I_3(\mathcal F_3(\bar S^2))$ construct the map
\begin{equation}\label{lamba}
\Lambda_x: \pi_2(K)/(i_1\pi_2(K_1)+i_2\pi_2(K_2)+i_3\pi_2(K_3))\to \frac{R_1\cap R_2\cap R_3}{[R_1,R_2\cap R_3][R_2,R_3\cap R_1][R_3,R_1\cap R_2]},
\end{equation}
where $F=\pi_1(K^1),$ $R_i=\ker\{F\to \pi_1(K_i)\},\ i=1,2,3,$ by setting
$$
\Lambda_x: (c_1,\dots, c_m)\mapsto [r_c, s_c].[R_1,R_2\cap R_3][R_2,R_3\cap R_1][R_3,R_1\cap R_2].
$$
Since $r_cs_ct_c=1$ in $F$, we have $[r_c,s_c]\in R_1\cap R_2\cap R_3$.

Let us show that the above map $\Lambda_x$ is well-defined. Let
$c'$ be an identity sequence equivalent to the sequence $c$.
Defining elements $r_{c'},s_{c'},t_{c'}$ as above, we have to show
that
\begin{equation}\label{equi}
[r_c,s_c]\equiv [r_{c'},s_{c'}]\mod [R_1,R_2\cap R_3][R_2,R_3\cap R_1][R_3,R_1\cap R_2].
\end{equation}
Since we above defined map $\Lambda_x$ is trivial for any sequence of the type (\ref{sequ3}) with conditions (\ref{ccond}), the equivalence
(\ref{equi}) is necessary and sufficient for the correctness of the map $\Lambda_x$.

First observe that if the sequences $c$ and $c'$ differ by the
Peiffer operations of the type (ii) or (iii), the equivalence
\ref{equi} holds. The only non-trivial Peiffer operations needed
to check are operations (iv) and (v). Since (v) is converse to
(iv), it is enough to prove the equivalence (\ref{equi}) for the
case $c'$ is obtained from $c$ by the single Peiffer operation of
the type (iv):
$$
c_{i}'=c_{i+1},\ c_{i+1}'=c_{i+1}^{-1}c_ic_{i+1},\ c_j'=c_j,\
j\neq i,i+1
$$
for some $1\leq i<m.$

The cases $i, i+1\in \{r_1,\dots, r_l\}$, $i, i+1\in \{s_1,\dots,
s_k\}$, $i,i+1\in \{t_1,\dots, t_h\}$ are trivial. In these cases
$r_c=r_{c'}, s_c=s_{c'},$ hence the needed equivalence
(\ref{equi}) follows. If $i+1\in \{r_1,\dots, r_l\},$ there is
also nothing to prove, since the definition of $r_c, s_c$ involves
the process of repeating of such operations. If $i\in \{t_1,\dots,
t_h\}$ or $i+1\in \{t_1,\dots, t_h\}$ then we clearly have
$[r_c,s_c]\equiv [r_{c'},s_{c'}] \mod [R_1,R_2\cap
R_3][R_2,R_3\cap R_1]$.

The only non-trivial case to consider is $i\in \{r_1,\dots,
r_l\},\ i+1\in \{s_1,\dots, s_k\}$. Clearly then,
$[r_{c'},s_{c'}]=[r_{c''},s_{c''}]$, where the sequence $c''$ is
obtained by applying again the operation (iv) to the sequence
$c':$
$$c_i''=c_{i+1}^{-1}c_ic_{i+1},\
c_{i+1}''=c_{i+1}^{-1}c_i^{-1}c_{i+1}c_ic_{i+1}.$$
Let $c_i=c_{r_j},$ $c_{i+1}=c_{s_e}$. Repeating the operation (iv), we can deform the
sequences $c$ and $c''$ to  the form
$$
r_2=r_1+1,\dots, r_{j-1}=r_{j-2}+1,\ r_{j+2}=r_{j+1}+1,\dots,
r_l=r_{l-1}+1
$$
without changing $[r_c,s_c]$ and $[r_{c''},s_{c''}]$. Now we can form the triple of words in $F$:
$$
\mathcal R_1'=\mathcal R_1\cup\{c_{r_1}\dots c_{r_{j-1}},
c_{j+1}\dots c_l\},\ \mathcal R_2'=\mathcal R_2,\ \mathcal
R_3'=\mathcal R_3.
$$
Clearly, this triple preserves the triple of normal subgroups
$R_1,R_2,R_3$ and we can consider the new identity sequences for the
triple of words $\mathcal R_1'\cup\mathcal R_2'\cup \mathcal R_3'$
formed by gluing the elements $c_{r_1},\dots, c_{r_{j-1}},$ and
$c_{r_{j+1}},\dots, c_{r_l}$: $$c'''=(*,\dots, *, c_{r_1}\dots
c_{r_{j-1}}, *,\dots, *, c_{r_{j+1}}\dots c_{r_l},*,\dots, *).$$
It is easy to see that
$$
[r_c,s_c]=[r_{c'''},s_{c'''}]
$$
in $F$. Hence, we can always assume that $l=3$, $c_{r_2}=c_i$ and
reduce arbitrary case to this one using the described procedure.
In these notations, we have sequences
\begin{align*}
& c=(*,\dots, *, c_{r_1},*, \dots, *, c_{r_2}, c_{s_e}, *,\dots,
*, c_{r_3},*,\dots, *),\\
& c''=(*,\dots, *, c_{r_1}, *, \dots, *, c_{r_2}^{c_{s_e}},
c_{s_e}^{c_{r_2}c_{s_e}},*,\dots,*, c_{r_3},*,\dots, *).
\end{align*}
We have the following:
\begin{align*}
& [r_c,s_c]=[c_{r_1}c_{r_2}c_{r_3}, S_1],\\
& [r_{c'''},s_{c'''}]=[c_{r_1}c_{r_2}^{c_{s_e}}c_{r_3}, S_2],
\end{align*}
where
\begin{align*}
& S_1=(\prod_{s_j<r_1}
c_{s_j}^{c_{r_1}c_{r_2}c_{r_3}})(\prod_{r_1<s_j<s_e}c_{s_j}^{c_{r_2}c_{r_3}})\cdot
c_{s_e}^{c_{r_3}}\cdot(\prod_
{s_e<s_j<r_3}c_{s_j}^{c_{r_3}})(\prod_{r_3<s_j}c_{s_j}),\\
& S_2=(\prod_{s_j<r_1}
c_{s_j}^{c_{r_1}c_{r_2}^{c_{s_e}}c_{r_3}})(\prod_{r_1<s_j<s_e}c_{s_j}^{c_{r_2}^{c_{s_e}}c_{r_3}})\cdot
c_{s_e}^{c_{r_2}c_{s_e}c_{r_3}}\cdot(\prod_
{s_e<s_j<r_3}c_{s_j}^{c_{r_3}})(\prod_{r_3<s_j}c_{s_j}).
\end{align*}
We then have
\begin{multline*}
[c_{r_1}c_{r_2}^{c_{s_e}}c_{r_3},
S_2]=c_{r_3}^{-1}c_{r_2}^{-1}[c_{r_2}^{-1},c_{s_e}]c_{r_1}^{-1}S_2^{-1}c_{r_1}c_{r_2}^{c_{s_e}}c_{r_3}S_2=\\
c_{r_3}^{-1}c_{r_2}^{-1}[c_{r_2}^{-1},c_{s_e}]c_{r_1}^{-1}S_2^{-1}c_{r_3}^{-1}c_{r_2}^{-c_{s_e}}c_{r_2}^{c_{s_e}}c_{r_3}c_{r_1}c_{r_2}^{c_{s_e}}c_{r_3}S_2
\equiv\\
c_{r_3}^{-1}c_{r_2}^{-1}c_{r_1}^{-1}S_2^{-1}c_{r_3}^{-1}c_{r_2}^{-c_{s_e}}[c_{r_2}^{-1},c_{s_e}]c_{r_2}^{c_{s_e}}c_{r_3}c_{r_1}c_{r_2}^{c_{s_e}}c_{r_3}S_2
\mod [R_3, R_1\cap R_2],
\end{multline*}
since $S_2^{-1}c_{r_3}^{-1}c_{r_2}^{-c_{s_e}}\in R_3,\
[c_{r_2}^{-1},c_{s_e}]\in R_1\cap R_2.$ Therefore,
$$
[c_{r_1}c_{r_2}^{c_{s_e}}c_{r_3},S_2]\equiv
c_{r_3}^{-1}c_{r_2}^{-1}c_{r_1}^{-1}S_2^{-1}c_{r_3}^{-1}c_{r_2}^{-c_{s_e}}c_{r_2}c_{r_3}c_{r_1}c_{r_2}^{c_{s_e}}c_{r_3}S_2\mod
[R_3,R_1\cap R_2].
$$
However,
$$
c_{r_1}c_{r_2}^{c_{s_e}}c_{r_3}S_2=c_{r_1}c_{r_2}c_{r_3}S_1\in R_3,
$$
therefore, $S_2=c_{r_3}^{-1}c_{r_2}^{-c_{s_e}}c_{r_2}c_{r_3}S_1$
and we have
$$
[c_{r_1}c_{r_2}^{c_{s_e}}c_{r_3},S_2]\equiv
[c_{r_1}c_{r_2}c_{r_3},S_1]\mod [R_3,R_1\cap R_2].
$$
Hence, we always have the needed equivalence (\ref{equi}) and we proved that the map $\Lambda_x$ is well-defined.
\end{proof}

For the generator $x\in \pi_3(S^2)$ denote by $\Lambda$ the composite map of the natural projection
$\pi_2(K)\to \pi_2(K)/i_1\pi_2(K_1)+i_2\pi_2(K_2)+i_3\pi_2(K_3)$ and the map $\Lambda_x$:
$$
\Lambda: \pi_2(K)\to \frac{R_1\cap R_2\cap R_3}{[R_1,R_2\cap R_3][R_2,R_3\cap R_1][R_3,R_1\cap R_2]}.
$$

\begin{prop}\label{sew1} Let $b\in i_{12}\pi_2(K_1\cup K_2)+i_{13}\pi_2(K_1\cup K_3)+i_{23}\pi_2(K_2\cup K_3)\subseteq \pi_2(K)$
where the maps $i_{12}, i_{13}, i_{23}$ are induced by the inclusions
$$
i_{12}: K_1\cup K_2\to K,\ i_{13}: K_1\cup K_3\to K,\ i_{23}: K_2\cup K_3\to K.
$$
We then have $\Lambda(a+b)=\Lambda(a)$ for every $a\in \pi_2(K)$.
\end{prop}
\begin{proof}
Let $a$ be an element from $\pi_2(K)$ presented by identity
sequence (\ref{div}) and the element $b$ be an element from
$i_{12}(K_1\cup K_2)\subseteq \pi_2(K)$ presented by the identity
sequence $(d_1,\dots, d_{l'}, e_{1},\dots, e_{k'})$ with $d_i\in
\mathcal R_1^F,\ e_i\in\mathcal R_2^F$. The element $a+b$ can be
presented by the following identity sequence
$$
c(a+b)=(c_{r_1},\dots, c_{r_l}, d_1,\dots, d_{l'}, c_{s_1}^{d_1\dots d_{l'}}, \dots, c_{s_k}^{d_1\dots d_{l'}}, e_1,\dots,e_{k'},f_1,\dots,f_{h'}),
$$
with $f_1,\dots, f_{h'}\in \mathcal R_3^F$. Denote
$a_1=c_{r_1}\dots c_{r_l}, a_2=d_1\dots d_{l'}, b_1=c_{s_1}\dots
c_{s_k},\ b_2=e_1\dots e_{k'}$. We then have
\begin{multline*}
[a_1a_2,b_1^{a_2}b_2]=a_2^{-1}a_1^{-1}b_2^{-1}a_2^{-1}b_1^{-1}a_2a_1b_1a_2b_2\equiv\\
a_1^{-1}b_2^{-1}a_2^{-1}b_1^{-1}a_1b_1a_2b_2\equiv a_1^{-1}b_1^{-1}a_1b_1\mod [R_3,R_1\cap R_2],
\end{multline*}
since $a_2\in R_1\cap R_2, a_1^{-1}b_2^{-1}a_2^{-1}b_1^{-1}\in R_3, a_2b_2=1$. Hence $\Lambda(a+b)=\Lambda(a)$.

In the case $b\in i_{13}\pi_2(K_1\cup K_3)+i_{23}\pi_2(K_2\cup K_3),$ we have obviously, that the elements which represent
$\Lambda(a+b)$ and $\Lambda(a)$ are equal modulo $[R_1,R_2\cap R_3][R_2, R_3\cap R_1]$ hence $\Lambda(a+b)=\Lambda(a)$.
\end{proof}

\vspace{.5cm}
The following example shows that the map $\Lambda$ is not always surjective.
\\

\noindent{\bf Example.} Let $F$ be a free group with generators $x_1,x_2$. Consider the following sets of words:
$$
\mathcal R_1=\{x_1\},\ \mathcal R_2=\{[x_1,x_2]\},\ \mathcal R_3=\{[x_1,x_2,x_1]\}.
$$
Denoting $R_1,R_2,R_3$ the normal closures of the sets $\mathcal R_1,\mathcal R_2, \mathcal R_3$ respectively, we have
$$
[R_1,R_2\cap R_3],[R_2,R_3\cap R_1], [R_3,R_1\cap R_2]\subseteq \gamma_4(F),
$$
where $\gamma_4(F)$ the 4-th lower central series term of $F$. However,
$$
[x_1,x_2,x_1]\in (R_1\cap R_2\cap R_3)\setminus \gamma_4(F),
$$
since $[x_1,x_2,x_1]$ is a basic commutator of length three in $F$. Suppose we have $$\Lambda(x)=[x_1,x_2,x_3].
[R_1,R_2\cap R_3][R_2,R_3\cap R_1][R_3,R_1\cap R_2]$$ for some element $x$ of the second homotopy module of the standard complex constructed
for the group presentation
$$
\langle x_1,x_2\ |\ x_1, [x_1,x_2],\ [x_1,x_2,x_1]\rangle.
$$
We then have \begin{equation}\label{poe}[x_1,x_2,x_1]\equiv
[r,s]\mod \gamma_4(F)\end{equation} for some $r\in R_1, s\in R_2$,
such that
\begin{equation}\label{po} rs\in R_3.
\end{equation}
However, the condition (\ref{po}) implies that $r\in \gamma_2(F),$ since $s\in \gamma_2(F)$. Therefore $[r,s]\in \gamma_4(F)$ and the equivalence
(\ref{poe}) is not possible. Hence, the map $\Lambda$ is
 not surjective.

\begin{theorem}\label{quadr}
The map $\Lambda$ is a homogenous quadratic map,
i.e. $$\Lambda(a,b)=\Lambda(a+b)-\Lambda(a)-\Lambda(b)$$ is bilinear and $\Lambda(x)=\Lambda(-x)$ for
any $a,b,x\in \pi_2(K)$.
\end{theorem}
\begin{proof}
For $x,y\in
\pi_2(K_{\langle X\ |\ \mathcal R_1\cup\mathcal R_2\cup\mathcal
R_3\rangle}),$ consider the cross-effect
$$
\Lambda(a,b)=\Lambda(a+b)-\Lambda(a)-\Lambda(b)\in \frac{R_1\cap R_2\cap R_3}{[R_1,R_2\cap
R_3][R_2,R_3\cap R_1][R_3,R_1\cap R_2]}.
$$
Represent elements $a,b$ by identity sequences:
$$
c(a)=(c_1,\dots, c_m),\ c(b)=(c_1',\dots, c_{m'}').
$$
Consider the corresponding divisions of the sequences $c(a)$ and $c(b)$:
\begin{align*}
& \{c_{r_1},\dots, c_{r_l}\}\cup\{c_{s_1},\dots,
c_{s_k}\}\cup\{c_{t_1},\dots, c_{t_n}\}=\{c_1,\dots, c_m\},\\
& \{c_{\bar r_1}',\dots, c_{\bar r_{l'}}'\}\cup\{c_{\bar
s_1}',\dots, c_{\bar s_{k'}}'\}\cup\{c_{\bar t_1}' ,\dots, c_{\bar
t_{n'}}'\}=\{c_1',\dots, c_m'\}
\end{align*}
with
$c_{r_i}, c_{\bar r_i}'\in \mathcal R_1^F, c_{s_i}, c_{\bar s_i}'\in \mathcal R_2^F,\ c_{t_i}, c_{\bar t_i}'\in \mathcal R_3^F.$
Consider then the induced division of the sequence
$c(a+b)=(c_1,\dots, c_m,c_1',\dots, c_{m'}'),$ which represents
the element $a+b\in \pi_2(K_{\langle X\ |\ \mathcal R_1\cup\mathcal
R_2\cup\mathcal R_3\rangle}):$
$$
\{c_{r_1},\dots, c_{r_l}, c_{\bar r_1}',\dots, c_{\bar
r_{l'}}'\}\cup \{c_{s_1},\dots, c_{s_k}, c_{\bar s_1}',\dots
c_{\bar s_{k'}}'\}\cup \{c_{t_1},\dots, c_{t_n}, c_{\bar
t_1}',\dots, c_{\bar t_{n'}}'\}.
$$
For the description of the functor $\Lambda(a,b)$, using the Peiffer operation
(iv) to the sequences $c(a)$ and $c(b)$, we can reduce the general
case to the case of $l=1, k=1, l'=1, k'=1$ with $r_1<s_1$, $\bar
r_1<\bar s_1$. Denote $x_1=c_{r_1}, y_1=c_{s_1}, x_2=c_{\bar
r_1}', y_2=c_{\bar s_1}'$.

We then have \begin{align*} &
\Lambda(a)=[x_1,y_1],\ \Lambda(b)=[x_2,y_2],\\
& \Lambda(a+b)=[x_1x_2, y_1^{x_2}y_2].
\end{align*}
We have
\begin{align*}
\Lambda(a+b)& = [x_1,y_2]^{x_2}[x_2,y_2][x_1,y_1^{x_2}]^{x_2y_2}[x_2,y_1^{x_2}]^{y_2}\\
& \equiv [x_1,y_2]^{x_2}[x_2,y_2][x_1,y_1^{x_2}][x_2,y_1]\mod [R_3,
R_1\cap R_2]\\
& \equiv
[x_1,y_2]^{x_2}[x_2,y_2]x_1^{-1}x_{2}^{-1}y_1^{-1}x_2x_1y_1\mod
[R_3,R_1\cap R_2]\\
& \equiv [x_1,y_2]^{x_2}[x_2,y_2][x_2,y_1]^{x_1}[x_1,y_1]\mod
[R_3,R_1\cap R_2].
\end{align*}
Since $x_1y_1, x_2y_2\in R_3,$
\begin{align*}
\Lambda(a,b)=\Lambda(a+b)-\Lambda(a)-\Lambda(b) & \equiv
[x_1,y_2]^{x_2}[x_2,y_1]^{x_1}\mod [R_3,R_1\cap R_2]\\ & \equiv
[x_1,y_2]^{y_2^{-1}}[x_2,y_1]^{y_1^{-1}}\mod [R_3,R_1\cap R_2]\\
& \equiv [y_2^{-1},x_1][y_1^{-1},x_2]\mod [R_3,R_1\cap R_2].
\end{align*}

Now let us show the linearity of the functor $\Lambda(*,*)$, i.e.
that
\begin{align}
& \Lambda(a+b,d)=\Lambda(a,c)+\Lambda(b,d),\label{f}\\
& \Lambda(a,b+d)=\Lambda(a,b)+\Lambda(a,d) \label{s}
\end{align}
for arbitrary elements $a,b,d\in \pi_2(K_{\langle X\ |\ \mathcal
R_1\cup\mathcal R_2\cup\mathcal R_3\rangle})$. Let $c(a), c(b)$ and
$c(d)$ be the identity sequences represented the elements $a, b$
and $d$ respectively. Again, without loss of generality we can
assume that these elements are represented by identity sequences with single element from each class $\mathcal R_i$.
Denote the correspondent pairs by $x_1,y_1\subset c(a)$ (the set-theoretical
inclusion means that $x_1$, $y_1$ are elements of the sequence
$c(a)$), $x_2,y_2\subset c(b),$ $x_3,y_3\subset c(d)$. In this
notation, modulo $[R_1,R_2\cap R_3][R_2, R_3\cap R_1][R_3,R_1\cap R_2]$, we have
\begin{align*}
\Lambda(a+b,d)& \equiv [y_3^{-1},x_1x_2][y_2^{-1}y_1^{-x_2},x_3]\\
& \equiv
[y_3^{-1},x_2][y_3^{-1},x_1]^{x_2}[y_2^{-1},x_3]^{y_1^{-x_2}}[y_1^{-x_2},x_3]\\
& \equiv
[y_3^{-1},x_2]x_2^{-1}y_3x_1^{-1}y_3x_1y_1x_2y_2x_3^{-1}y_2^{-1}x_2^{-1}y_1^{-1}x_2x_3\\
& \equiv
[y_3^{-1},x_2]x_2^{-1}y_3x_1^{-1}y_3x_1y_1(x_2y_2x_3^{-1}y_2^{-1}x_2^{-1}x_3)x_3^{-1}y_1^{-1}x_2x_3\\
& \equiv
[y_3^{-1},x_2]x_2^{-1}(x_2y_2x_3^{-1}y_2^{-1}x_2^{-1}x_3)y_3x_1^{-1}y_3x_1y_1x_3^{-1}y_1^{-1}x_2x_3\\
& \equiv
[y_3^{-1},x_2][y_2^{-1},x_3]x_3^{-1}x_2^{-1}x_3y_3x_1^{-1}y_3x_1y_1x_3^{-1}y_1^{-1}x_2x_3\\
& \equiv
[y_3^{-1},x_2][y_2^{-1},x_3]x_3^{-1}x_2^{-1}x_3[y_3^{-1},x_1][y_1^{-1},x_3]x_3^{-1}x_2x_3\\
& \equiv
[y_3^{-1},x_2][y_2^{-1},x_3][y_3^{-1},x_1][y_1^{-1},x_3]\\
& \equiv \Lambda(a,d)+\Lambda(b,d),
\end{align*}
since $[y_3^{-1},x_1][y_1^{-1},x_3]\in R_2\cap R_3$ and (\ref{f})
follows. The equality (\ref{s}) can be proved analogically.

Now let us prove that $\Lambda(-x)=\Lambda(x).$ Clearly, we can assume that
our identity sequence representing the element $x\in \pi_2(K)$ has the form
$$
(r_1,s_1,t_1)
$$
with $r_1\in \mathcal R_1,s_1\in \mathcal R_2,t_1\in \mathcal R_3$. The
inverse sequence, which represents the element $-x$ has the form
$$
(t_1^{-1},s_1^{-1},r_1^{-1}).
$$
We have
$$
\Lambda(-x)=[r_1^{-1},
s_1^{-r_1^{-1}}]=[s_1^{-1},r_1]=[r_1,s_1]^{s_1^{-1}}\equiv
[r_1,s_1]\equiv \Lambda(x)\mod [R_2, R_3\cap R_1].
$$
\end{proof}

\begin{theorem}\label{thmain}
The function $\Lambda$ induces the homomorphism of
$F/R_1R_2R_3$-modules
$$
\bar \Lambda: \pi_3(K)\to \frac{R_1\cap R_2\cap R_3}{[R_1,R_2\cap R_3][R_2,R_3\cap
R_1][R_3,R_1\cap R_2]}.
$$
\end{theorem}
\begin{proof}
Let $x\in \pi_2(K_{\langle X\ |\ \mathcal R_1\cup \mathcal R_2\cup
\mathcal R_3\rangle})$. Present $x$ by the sequence
$$
c(x)=(c_1,\dots, c_m).
$$
For a given element $f\in \pi_1(K),$ present this element as a
coset $f=w.R_1R_2R_3$ for some element $w\in F$. The element
$f\circ x\in \pi_2(K_{\langle X\ |\ \mathcal R_1\cup \mathcal
R_2\cup \mathcal R_3\rangle})$ can be presented by sequence
$$
c(x)^w=(c_1^w,\dots, c_m^w).
$$
It follows directly from the definition of $\Lambda(x),$ that
$$
\Lambda(f\circ x)\equiv \Lambda(x)^w\mod [R_1,R_2\cap R_3][R_2,R_3\cap R_1][R_3,R_1\cap R_2].
$$
Since $\pi_3(K)=\Gamma \pi_2(K)$, we have the needed homomorphism
of $F/R_1R_2R_3$-modules due to Theorem \ref{quadr}.
\end{proof}

\noindent{\bf Example.} For two-dimensional sphere $S^2$, clearly,
$\Lambda$ defines the isomorphism (\ref{wujie}):
$$
\bar\Lambda: \pi_3(S^2)\to I_3(\mathcal F_3(\bar S_3))
$$
with $\bar S_3\in\mathcal K_3$ defined in (\ref{sbar}).\\

\noindent{\bf Example.} Consider a group presentation
$$
\mathcal P=\langle x_1,\dots, x_k\ |\ r_1,\dots, r_l\rangle
$$
of a group $G$. Let $\mathcal P'$ be another presentation of $G$ with $k+2l$ generators and $3l$ relators
given by
$$
\mathcal P'=\langle x_1,\dots, x_k, y_1,\dots, y_l, z_1,\dots, z_l\ |\ y_1,\dots, y_l, z_1y_1^{-1},\dots, z_ly_l^{-1}, z_1^{-1}r_1,\dots, z_l^{-1}r_l\rangle
$$
The standard complex $K_{\mathcal P'}$ is the union $K_1\cup
K_2\cup K_3$, where $K_1,K_2,K_3$ are standard complexes of the
following presentations
\begin{align*}
& \langle x_1,\dots, x_k, y_1,\dots, y_l,z_1,\dots, z_l\ |\ y_1,\dots, y_l\rangle,\\
& \langle x_1,\dots, x_k, y_1,\dots, y_l,z_1,\dots, z_l\ |\ z_{1}y_1^{-1},\dots, z_ly_l^{-1}\rangle,\\
& \langle x_1,\dots, x_k, y_1,\dots, y_l,z_1,\dots, z_l\ |\ z_1^{-1}r_1,\dots, z_l^{-1}r_l\rangle
\end{align*}
respectively. Denoting $\bar K=(K_{\mathcal P'},K_1,K_2,K_3)\in \mathcal K_3$, we have the following isomorphism of $G$-modules:
$$
\pi_3(K_{\mathcal P})\simeq \pi_3(K_{\mathcal P'})\simeq I_3(\mathcal F_3(\bar K)).
$$
This isomorphism follows directly from the description of Kan's
loop construction $GK_{\mathcal P}$ and the fact that for a
simplicial group $G_*$ with $G_2$ generated by degeneracy
elements, one has $\pi_2(G_*)\simeq I_3(G_2, ker(d_0), ker(d_1),
ker(d_2))$ (see, for example, \cite{MP}). \vspace{.5cm}

\section{Application to group homology}
\vspace{.5cm} \subsection{} For a given element
$(F;R_1,\dots,R_n)\in \mathcal R_n$ it follows from \cite{BE} and
\cite{DIP} that under certain conditions on $(F;R_1,\dots, R_n)$,
the group homologies, or more generally, the derived functors of
the lower central quotients, can be obtained with the help of the
groups $I_n(F;R_1,\dots, R_n)_F$, i.e. the $F$-coinvariant part of
$I_n(F;R_1,\dots,R_n)$.

For every $(F;R_1,R_2)\in \mathcal R_2$, it is well-known that
there exists a canonical map
\begin{equation}\label{ert}
H_3(G)\to I_2(F;R_1,R_2)=\frac{R_1\cap R_2}{[R_1,R_2][F,R_1\cap
R_2]}
\end{equation}
which is a part of a long exact sequence of homology groups, see
\cite{BG} and \cite{DEG}. This map can be easily obtained from the
Gutierrez-Ratcliffe map (\ref{exa1}). For that we consider
arbitrary $(K,K_1,K_2)\in \mathcal K_2$ with $F=\pi_1(K^1),\
R_i=ker\{\pi_1(K^1)\to \pi_1(K_i)\},\ i=1,2$. The chain complex
\begin{equation}\label{cell}0\to \pi_2(K)\to C_2(\tilde K)\to C_1(\tilde K)\to \mathbb
Z[F/R_1R_2]\to \mathbb Z\to 0\end{equation} of the universal cover
$\tilde K$ of $K$ can be viewed as a complex of free
$F/R_1R_2$-modules. Applying the group homology functor
$H_*(F/R_1R_2,-)$ to (\ref{cell}), we obtain the natural maps
$$
\partial_n: H_n(F/R_1R_2)\to H_{n-3}(F/R_1R_2,\pi_2(K)),\ n\geq 3,
$$
which are isomorphisms for $n\geq 4$. The map (\ref{ert}) is the
composition of $\partial_3$ and the $F$-coinvariant map from
(\ref{exa1}).

\subsection{} Recall the definitions of certain quadratic functors in
the category of abelian groups. Let $A$ be an abelian group.
Define the symmetric tensor square
$$
SP^2(A)=A\otimes A/\{a\otimes b-b\otimes a,\ a,b\in A\},
$$
and the augmentation power functor
$$
P_2(A)=\Delta(A)/\Delta^3(A),\ \Delta(A)=ker\{\mathbb Z[A]\to
\mathbb Z\}.
$$
It is well-known (see, for example, \cite{Baues1}) that for a free
abelian group $A$, there are the following short exact sequences
of abelian groups:
\begin{align}
& 0\to SP^2(A)\to P_2(A)\to A\to 0 \label{o1}\\
& 0 \to SP^2(A)\to \Gamma(A)\to A\otimes \mathbb Z_2\to 0
\label{o2}.
\end{align}
Now let $A$ be a $G$-module. The $G$-action can be naturally
extended to the abelian groups $SP^2(A),P_2(A),\Gamma(A),A\otimes
\mathbb Z_2$, thus allowing to consider the sequences (\ref{o1})
and (\ref{o2}) as sequences of $G$-modules. Applying homology
functor $H_*(G,-)$ to (\ref{o1}) and (\ref{o2}) we obtain long
exact sequences
\begin{align*}
& \dots \to H_1(G, P_2(G))\to H_1(G,A)\to H_0(G,SP^2(A))\to H_0(G,
P_2(A))\to \dots\\
& \dots \to H_1(G,\Gamma(A))\to H_1(G, A\otimes \mathbb Z_2)\to
H_0(G,SP^2(A))\to H_0(G,\Gamma(A))\to \dots
\end{align*}
\subsection{}
It is natural to ask about applications of the map $\bar\Lambda$
constructed in Theorem \ref{thmain}. To this end, we consider $(K;
K_1,K_2,K_3)\in \mathcal K_3$ with $$\pi_1(K^1)=F,\
R_i=ker\{\pi_1(K^1)\to \pi_1(K_i)\},\ i=1,2,3.$$ Let us denote
$G=F/R_1R_2R_3$ and let us define the map
$$
\Psi_4: H_4(G)\to I_3(F;R_1,R_2,R_3)_F=\frac{R_1\cap R_2\cap
R_3}{[R_1,R_2\cap R_3][R_2,R_3\cap R_1][R_3,R_1\cap R_2][F,R_1\cap
R_2\cap R_3]}
$$
as a composite map in the following diagram with exact rows and
columns:
$$
\xyma{ & H_4(G) \ar@{=}[d] & H_1(G,\pi_2(K)\otimes \mathbb Z_2) \ar@{->}[d]\\
H_1(G, P_2(\pi_2(K))) \ar@{->}[r] & H_1(G,\pi_2(K)) \ar@{->}[r] &
H_0(G,SP^2(\pi_2(K)))
\ar@{->}[d]\\
& H_0(G,\pi_3(K)) \ar@{=}[r] \ar@{->}[d]_{H_0(F,\bar\Lambda)} & H_0(G,\Gamma\pi_2(K)) \ar@{->>}[d]\\
& I_3(F;R_1,R_2,R_3)_F & H_0(G,\pi_2(K)\otimes \mathbb Z_2)}
$$
\noindent{\bf Remark.} Proposition \ref{sew1} implies that the
natural composition map
$$
H_4(F/R_1R_2)\oplus H_4(F/R_2R_3)\oplus H_4(F/R_1R_3)\to
H_4(G)\buildrel{\Psi_4}\over\rightarrow I_3(F; R_1,R_2,R_3)_F
$$
is the zero map.

\end{document}